\documentclass[11pt, a4paper, reqno]{amsart}

\usepackage{amsfonts, amsthm, amsmath, amscd, amssymb}

\newcommand\QQ{\mathbb{Q}}

\newcommand\NN{\mathbb{N}}
\newcommand\FF{\mathbb{F}}
\newcommand\RR{\mathbb{R}}

\newcommand\ZZ{\mathbb{Z}}

\newcommand\PP{\mathbb{P}}
\newcommand\x{\mathbf{x}}

\newcommand\Ptwo{{\PP^2}}

\DeclareMathOperator{\hcf}{gcd}

\renewcommand{\le}{\leqslant}
\renewcommand{\ge}{\geqslant}
\renewcommand{\leq}{\leqslant}
\renewcommand{\geq}{\geqslant}

\newcommand\la{\lambda}
\newcommand\al{\alpha}

\newcommand\sfl{\mathsf{\Lambda}}
\newcommand\sfm{\mathsf{M}}

\newtheorem{thm}{Theorem}
\newtheorem*{thm*}{Theorem}
\newtheorem*{lemma}{Lemma}

\theoremstyle{definition}
\newtheorem*{ack}{Acknowledgement}

\begin{document}

\title{Solubility of Fermat equations}

\author{T.D. Browning}
\author{R. Dietmann}

\address{Institut f\"ur Algebra und Zahlentheorie,
Lehrstuhl f\"ur Zahlentheorie,
Pfaffenwaldring 57,
D-70569 Stuttgart}
\email{dietmarr@mathematik.uni-stuttgart.de}

\address{School of Mathematics,  
University of Bristol, Bristol BS8 1TW}
\email{t.d.browning@bristol.ac.uk}


\date{\today}

\begin{abstract}
The arithmetic of the equation 
$a_1 x_1^d + a_2 x_2^d + a_3 x_3^d = 0$ is considered for $d \geq 2$, with the
outcome that the set of coefficients for which the equation admits a
non-zero integer solution is shown to have density zero.  
\end{abstract}

\maketitle

\section{Introduction}

Let $d\in\NN$.  The purpose of this short note is to 
discuss the locus of rational points $C_d(\QQ)$ on the Fermat curves
\begin{equation}\label{eq:fermat}
C_d:\quad   a_1 x_1^d + a_2 x_2^d + a_3 x_3^d = 0
\end{equation}
in $\Ptwo$, for given $\mathbf{a}=(a_1,a_2,a_3)\in \ZZ^3$. 
Our main goal is to show that a random such curve does not possess a
rational point for $d>1$. However, we will also discuss the number of 
solutions when  $C_d(\QQ)\neq
\emptyset$, and briefly consider the analogous problem in higher dimension.

Throughout our work we will allow $i,j,k$ to denote distinct elements
from the set $\{1,2,3\}$.
For any $H\geq 1$, let $N_d(H)$ denote the 
number of $\mathbf{a}\in \ZZ^3$ with
$|a_i| \le H$,  for which $C_d(\QQ)\neq
\emptyset$.  The following is our main result.

\begin{thm}\label{t:main}
We have 
\[
  N_d(H) \ll \frac{H^3}{(\log H)^{\psi(d)}},
\]
where if $\phi$ denotes Euler's totient function then 
$$
\psi(d):=\frac{3}{\phi(d)}\Big(1-\frac{1}{d}\Big).
$$
\end{thm}


All of the implied constants in our work are allowed to depend at
most upon $d$.  In the case $d=2$ of conics, Theorem \ref{t:main} yields
$$
  N_2(H) \ll \frac{H^3}{(\log H)^{\frac{3}{2}}}.
$$
This retrieves an earlier result of Serre \cite{serre}. 
In fact Guo \cite{g} has established an asymptotic formula for the
corresponding quantity in which the coefficients $a_i$ are restricted
to be odd, with $a_1a_2a_3$ square-free.
In the case of additive cubics our result implies that 
$$
  N_3(H) \ll \frac{H^3}{\log H}.
$$
This provides a partial answer to a question raised by Poonen and
Voloch \cite{p-v}:  does a random cubic curve in $\Ptwo$ that is
defined over $\QQ$ possess a
$\QQ$-rational point?
The proof of Theorem \ref{t:main} will be established in \S
\ref{s:large} using the large sieve
inequality, as pioneered by Serre \cite{serre} in the case $d=2$.

It is natural to ask what happens for additive equations in more than
three variables. 
For given $d\in \NN$ and $H\geq 1$, let $M_d(H)$ denote the number of
$\mathbf{a} \in \mathbb{Z}^4$ with $|a_i| \le H$, for which the equation 
$$
a_1 x_1^d + a_2 x_2^d +a_3 x_3^d + a_4 x_4^d = 0
$$
is everywhere locally soluble. 
The following inequality will be proved in \S \ref{s:lift}.

\begin{thm}
\label{quat}
We have 
$$
M_d(H) \gg H^4.
$$
\end{thm}

Theorem \ref{quat} provides an additive analogue of a result due to 
Poonen and Voloch \cite[Theorem 3.6]{p-v}. The latter establishes that 
a positive proportion of all hypersurfaces in $\mathbb{P}^{n-1}$ of
degree $d$ that are defined over $\QQ$ are everywhere locally soluble,
provided that $n-1,d\geq 2$ and $(n,d)\neq (3,2)$.

Returning to the setting of ternary forms, 
let us now consider the problem of describing $C_{d}(\QQ)$ when it is
non-empty. When $d$ is sufficiently large, it has been
conjectured by Granville \cite{granville} on the basis of
a generalised version of the $abc$-conjecture, that the curve
\eqref{eq:fermat} never has any non-trivial rational points.  
While Faltings' proof of the Mordell conjecture ensures that there are
only finitely many rational points for each $d\geq 4$, 
it is notoriously difficult to achieve an effective bound for the total number of solutions in
terms of the coefficients $a_1,a_2,a_3$.  The following result deals
with the much simpler scenario in which one restricts attention to
rational points of bounded height on the curve. 

\begin{thm} \label{upper}
Let $d\geq 2$ and let $\mathbf{a} \in \ZZ^3$ 
have pairwise coprime non-zero components. Then we have 
$$
\#\{x \in C_d(\QQ) : H(x)\leq B\} \ll 
\Big(1+\frac{B^{\frac{3}{2}}}{|a_1a_2a_3|^{\frac{1}{d}}}\Big)
d^{\omega(a_1a_2a_3)}, 
$$
where $H:\Ptwo(\QQ)\rightarrow \RR_{\geq 0}$ is the usual exponential height function.
\end{thm}

Theorem \ref{upper} will be established in \S \ref{s:geom}. 
It generalises a result due to Heath-Brown \cite[Theorem 2]{hb} which
deals with the case $d=2$. Theorem \ref{upper} is 
susceptible to improvement in a number of obvious directions. 
Firstly it would be easy to extend this result to counting rational
points whose coordinates are constrained to lie in
lopsided boxes, rather than in a cube.
Secondly, at the expense of weakening the dependence of the estimate
on $a_1,a_2,a_3$, the exponent of $B$ can be improved substantially. 
We will not pursue either of these lines of enquiry here, however.
Theorem \ref{upper} can be used to provide some simple-minded
evidence for the expected paucity of rational points on
\eqref{eq:fermat}.  
Thus when $a_1,a_2,a_3$ are arbitrary integers, it follows from this result that there are only
$O(d^{\omega(a_1a_2a_3)})$ points in $C_d(\QQ)$ with height at most $|a_1a_2a_3|^{2/(3d)}$.
Thus if there are many rational points then their height must be large
compared to the height of the defining form.

\begin{ack}
While working on this paper the first author was supported by EPSRC
grant number \texttt{EP/E053262/1}.
\end{ack}

\section{The large sieve}\label{s:large}


For any $H\geq 1$ we let $N_d^*(H)$ be defined as for $N_d(H)$ but
with the extra hypotheses that $a_1a_2a_3\neq 0$ and $\gcd(a_1,a_2,a_3)=1$.
We will say that an integer $a\in \NN$ is $d$-free if $\nu<d$ whenever
$p^\nu\mid a$. For $H_i\geq 1$ we let $N_d^{**}(\mathbf{H})$ denote
the number of 
$d$-free triples $\mathbf{a}\in \ZZ^3$ with
$0<|a_i| \le H_i$ and $\gcd(a_1,a_2,a_3)=1$, such that
$C_d(\QQ)\neq
\emptyset$.  We will use the large sieve inequality to show that
\begin{equation}
  \label{eq:show}
  N_d^{**}(\mathbf{H}) \ll \frac{\prod_{i=1}^3(H_i+z^2)}{(\log z)^{\psi(d)}},
\end{equation}
for any $z \geq 1$, where $\psi(d)$ is given in the statement of
Theorem \ref{t:main}.

Let us begin by seeing how this suffices for the
statement of Theorem \ref{t:main}. Now it is clear that 
$$
N_d(H)=\sum_{k\leq H} N_d^*(k^{-1}H)+O(H^2),
$$
whence it certainly suffices to establish the theorem for 
$N_d^*(H)$ in place of $N_d(H)$. We now write $a_i=u_iv_i^d$ in the
definition of $N_d^*(H)$, with each $u_i$ being $d$-free and $v_i>0$. It 
follows that 
$$
N_d^*(H)=\sum_{v_i\leq H^{\frac{1}{d}}} N_d^{**}(\mathbf{H}),
$$
where $\mathbf{H}$ has components $H_i:=H/v_i^d$. 
Let us break the summation over $\mathbf{v}\in \NN^3$ into 
two set $S_1(H)\cup S_2(H)$, where $S_1(H)$ denotes the set of vectors 
for which one of the components is bigger than $(\log H)^{3/(d\phi(d))}$,  and
$S_2(H)$ denotes the remainder.   It is trivial to see that 
\begin{align*}
\sum_{\mathbf{v} \in S_1(H)} N_d^{**}(\mathbf{H})
\leq 
\sum_{\mathbf{v} \in S_1(H)} \#\{\mathbf{u}\in\ZZ^3: |u_i|\leq H_i\}
&\ll
H^3\sum_{\mathbf{v} \in S_1(H)} \frac{1}{v_1^dv_2^dv_3^d}\\
&\ll 
\frac{H^3}{(\log H)^{\psi(d)}},
\end{align*}
which is satisfactory. 
Turning to the contribution from 
the set $S_2(H)$, we deduce from \eqref{eq:show} that 
\begin{align*}
\sum_{\mathbf{v} \in S_2(H)} N_d^{**}(\mathbf{H})
&\ll  \sum_{\mathbf{v} \in S_2(H)} 
\frac{\prod_{i=1}^3(H_i+z^2)}{(\log z)^{\psi(d)}},
\end{align*}
for any $z\geq 1$. Since $v_i\leq (\log H)^{3/(d\phi(d))}$, so it follows that
$$
H_i\geq \frac{H}{v_1^dv_2^dv_3^d} \gg 1,
$$ 
for $H\gg 1$.
Taking $z=H^{1/2}(v_1v_2v_3)^{-d/2}$ therefore yields 
\begin{align*}
\sum_{\mathbf{v} \in S_2(H)} N_d^{**}(\mathbf{H})
&\ll  \frac{H^3}{(\log H)^{\psi(d)}}\sum_{\mathbf{v} \in S_2(H)} 
\frac{1}{v_1^dv_2^dv_3^d}
\ll 
\frac{H^3}{(\log H)^{\psi(d)}}.
\end{align*}
This completes the deduction of Theorem \ref{t:main}, subject to \eqref{eq:show}.

We now proceed with the proof of \eqref{eq:show}. 
Let $p>2$ and let $R_d(p)$ denote the number of 
$a\in \FF_p^*$ for which there is a solution of the congruence
$$
x^d\equiv a \pmod{p}.
$$
It is an easy exercise in elementary number theory (see \cite[\S 4.2]{i-r},
for example) to show that 
\begin{equation}
  \label{eq:Md}
R_d(p)=\frac{p-1}{\gcd(d,p-1)}.
\end{equation}
We will be interested in the set of $\mathbf{a}\in \FF_p^3$ which arise
as images of the points counted by 
$N_d^{**}(\mathbf{H})$. We denote the cardinality of such
$\mathbf{a}\in \FF_p^3$ by $p^3-\tau(p)$, where $\tau(p)$ denotes the
number of vectors in $\FF_p^3$ that are excluded.

We seek a good lower bound for $\tau(p)$, still under the assumption
that $p>2$. Let $\mathbf{a}\in \FF_p^3$ be such that $p\mid
a_i$ and $p\nmid a_ja_k$. Then for 
fixed $a_j\in \FF_p^*$ there are exactly $R_d(p)$ values
of $a_k$ for which the congruence 
\[
  a_j x_j^d + a_k x_k^d \equiv 0 \pmod p
\]
has solutions with $p\nmid x_jx_k$.  For the remaining vectors 
with $p\mid a_i$ and $p\nmid a_ja_k$, the only solution to the above
congruence has $p\mid \gcd(x_j,x_k)$. But then the condition
$C_d(\QQ)\neq \emptyset$ in $N_d^{**}(\mathbf{H})$ implies that 
$$
a_ix_i^d\equiv 0\pmod{p^d},
$$
for some $x_i\in \ZZ$ which is coprime to $p$. 
This is impossible since $a_i$ is $d$-free.
Employing \eqref{eq:Md} this therefore 
establishes that 
\begin{align*}
\tau(p)\geq \sum_{i} 
(p-1)^2\Big(1-\frac{1}{\gcd(d,p-1)}\Big)
&= 3(p-1)^2\Big(1-\frac{1}{\gcd(d,p-1)}\Big)\\
&\geq  3(p-1)^2\Big(1-\frac{1}{k}\Big)\\
&=  3p^2\Big(1-\frac{1}{k}\Big)+O(p),
\end{align*}
for any $k\mid \gcd(d,p-1)$.

We are now ready for our application of the large sieve inequality in
dimension three. Let $z,H_i\geq 1$. It easily follows from 
the arguments of \cite{large} that
\begin{equation}
  \label{eq:**}
N_d^{**}(\mathbf{H})\ll 
\frac{\prod_{i=1}^3(H_i+z^2)}{G(z)},
\end{equation}
where
$$
G(z):=\sum_{\substack{n \leq z\\ 2\nmid n}}|\mu(n)|\prod_{p \mid n} \frac{\tau(p)}{p^3-\tau(p)}.
$$
For any $k\in \NN$, let $\mathcal{P}_k$ denote the set of primes
congruent to $1$ modulo $k$, and let 
$g_k$ be the non-negative multiplicative arithmetic function
$$
g_k(n):=\frac{|\mu(n)|\big(3(1-\frac{1}{k})\big)^{\omega(n)}}{n}.
$$
Then we have 
\begin{align*}
G(z)
\geq 
\sum_{\substack{k\mid d\\ k>1}}
\sum_{\substack{n \leq z\\ p\mid n \Rightarrow p \in \mathcal{P}_k}}
|\mu(n)|\prod_{p \mid n} \left(\frac{3}{p}\Big(1-\frac{1}{k}\Big)
  +O\Big(\frac{1}{p^2}\Big)\right)
&\gg
\sum_{\substack{k\mid d\\ k>1}}
\sum_{\substack{n \leq z\\ p\mid n \Rightarrow p \in
    \mathcal{P}_k}}g_k(n)\\
&\geq
\sum_{\substack{n \leq z\\ p\mid n \Rightarrow p \in
    \mathcal{P}_d}}g_d(n).
\end{align*}
It will be convenient to set $\gamma:=3(1-1/d)$. 
Now it is easy to see that 
$$
\sum_{n \leq z}g_d(n) \gg (\log z)^{\gamma},
$$
and furthermore,
\begin{align*}
\sum_{\substack{n \leq z\\ p\mid n \Rightarrow p \not\in
    \mathcal{P}_d}}g_d(n) \leq \exp \Big(
\sum_{ \substack{p \leq z\\ p \not\in
    \mathcal{P}_d}}\frac{\gamma}{p}\Big)
\ll (\log z)^{\gamma(1-\frac{1}{\phi(d)})},
\end{align*}
by Dirichlet's theorem on primes in arithmetic progression. Hence it
follows that
$$
G(z)\gg \sum_{\substack{n \leq z\\ p\mid n \Rightarrow p \in
    \mathcal{P}_d}}g_d(n)  \gg 
(\log z)^{\gamma} \Big(\sum_{\substack{n \leq z\\ p\mid n \Rightarrow p \not\in
    \mathcal{P}_d}}g_d(n)\Big)^{-1}\gg (\log z)^{\frac{\gamma}{\phi(d)}},
$$
in \eqref{eq:**}.
On noting that $\gamma/\phi(d)=\psi(d)$, this therefore 
completes the proof of \eqref{eq:show}.

\section{A lifting argument}\label{s:lift}

In this section we prove  Theorem \ref{quat}, which will be achieved
via a simple lifting argument. As is well known, there exists a
constant $p_0=p_0(d)>d$ 
such that for  primes $p \ge p_0$ we will have $p \nmid d$, and
furthermore, every congruence of the form
\[
  b_1 x_1^d + \cdots + b_3 x_3^d \equiv 0 \pmod p
\]
will have a non-trivial solution when $p \nmid b_1 b_2 b_3$.
Any such solution can be lifted to a non-trivial solution in 
$\mathbb{Q}_p^3$. For each prime $p < p_0$ there exists a power of
$p$, which we denote by $q_p \in
\mathbb{N}$, and a residue class $\mathbf{a}_p$ modulo $q_p$, such that 
the congruence
\[
  a_1 x_1^d + \cdots + a_4 x_4^d \equiv 0 \pmod {q_p}
\]
has a solution which can be lifted to a non-trivial solution
in $\mathbb{Q}_p^4$, if 
$$
\mathbf{a} \equiv \mathbf{a}_p \pmod {q_p}.
$$ 
Let
$Q := \prod_{p < p_0} q_p.$
By the Chinese remainder theorem there exists a residue class
$\mathbf{c} \bmod Q$ such that $\mathbf{a} \equiv \mathbf{c} \bmod{Q}$ implies
that $\mathbf{a} \equiv \mathbf{a}_p \bmod {q_p}$ for $p < p_0$.
It follows that $M_d(H) \ge M_d^*(H)$,  where
$M_d^*(H)$ denotes the number of $\mathbf{a}\in \mathbb{Z}^4$
such that $|a_i| \le H$, with 
$
\mathbf{a} \equiv \mathbf{c} \bmod{Q}
$, 
and furthermore, any $p \ge p_0$ 
divides at most one of the coefficients $a_i$.
Assuming that $p_0$ is taken to be sufficiently large we deduce that $M_d^*(H) \gg H^4$, which therefore 
concludes the proof of Theorem~\ref{quat}.

\section{Geometry of numbers}\label{s:geom}

In establishing Theorem \ref{upper} we will need to establish the existence of a small
number of lattices, each of reasonably large determinant, that can be
used to cover the integer solutions to the equation in \eqref{eq:fermat}.
This is provided by the following result.

\begin{lemma}
Let $d \geq 2$ and let $\mathbf{a}\in \ZZ^3$ such that $a:=|a_1a_2a_3|\neq 0$ and 
$\hcf(a_i,a_j)=1$.
Let $\x\in\ZZ^3$ be a solution of \eqref{eq:fermat}. 
Then there exist lattices $\sfl_1,\ldots, \sfl_J \subseteq \ZZ^3$ such
that 
\begin{enumerate}
\item[(i)] $J\leq 2d^{\omega(a)}$,
\item[(ii)] $\dim \sfl_j=3$ and $\det \sfl_j \gg a^{2/d}$, for each $j\leq J$,
\item[(iii)] $\x\in\bigcup_{j\leq J}\sfl_j$.
\end{enumerate}
\end{lemma}

\begin{proof}
Let $p\mid a$, and write $\al_i=v_p(a_i)$.  
Since the integers $a_1,a_2,a_3$ are pairwise coprime, we may assume
without loss of generality that $\al_1=\al_2=0$ and $\al_3\geq 1.$
 For any prime $p$ let
$\delta=v_p(d)$ be the $p$-adic order of $d$, and write
$$
  \gamma:=\begin{cases}
\delta+1, & \mbox{if $p>2$ or $\delta=0$},\\
\delta+2, & \mbox{if $p=2$ and $\delta\geq 1$.}
\end{cases}
$$
Let $\x\in\ZZ^3$ be such that \eqref{eq:fermat} holds. We claim that
there exist sublattices $\sfm_1,\ldots,
\sfm_K\subseteq\ZZ^3$ with  $K\leq 2^{\gamma-\delta-1} d$, such that
$$
\dim \sfm_k=3, \quad \det \sfm_k \geq p^{\frac{2\al_3}{d}-\gamma+1},
$$ 
for each $k\leq K$, and $\x\in\bigcup_{k\leq K}\sfm_k$.
The Chinese remainder theorem will then produce at most $2d^{\omega(a)}$ integer
sublattices overall, each of dimension $3$ and determinant
\begin{align*}
\geq \prod_{p\mid a_1a_2a_3} p^{\frac{2v_p(a)}{d}-\gamma+1} 
= a^{\frac{2}{d}} \prod_{p\mid a_1a_2a_3} p^{-\gamma+1}
\gg a^{\frac{2}{d}}.
\end{align*}
This completes the proof of
the lemma subject to the construction of the lattices $\sfm_1,\ldots,\sfm_K$.

Turning to the claim, 
let $\x\in\ZZ^3$ be such that \eqref{eq:fermat} holds. 
Let us write $x_i=p^{\xi_i}x_i'$, for
$i=1,2$, with $p\nmid x_1'x_2'$ and $\xi_1\leq \xi_2$, say. Then 
we deduce that  
$$
a_1p^{d\xi_1}{x_1'}^d+a_2p^{d\xi_2}{x_2'}^d\equiv 0 \pmod{p^{\al_3}}.
$$
There are now $3$ possibilities to consider: either $\al_3\leq
d\xi_1$, or $d\xi_1<\al_3\leq d\xi_2$, or $d\xi_2<\al_3$. 
The second case is plainly impossible.  In the first case we may conclude that $\x$
belongs to the set of $\x\in\ZZ^3$ such that $p^{\lceil \al_3/d\rceil
}$ divides $x_1$ and $x_2$. This defines an integer lattice 
of dimension $3$ and determinant $\geq p^{2\al_3/d}$.  Thus we may
take $K=1$ in this case.  

Finally, in the third case, we must have $\xi_1=\xi_2=\xi$, say. But
then it follows that
$$
a_1{x_1'}^d+a_2{x_2'}^d\equiv 0 \pmod{p^{\al_3-d\xi}}.
$$
Suppose first that $\al_3-d\xi<\gamma$.  Then we have
$$
p^{\frac{2\al_3}{d}}\leq p^{2\xi}p^{\frac{2(\gamma-1)}{d}}\leq
p^{2\xi+\gamma-1}. 
$$
Since $\x$ lies on the lattice of determinant $p^{2\xi}$ that is
determined by the conditions $p^\xi \mid x_1$ and $p^\xi \mid x_2$, we
may clearly take $K=1$ in this case also.
Suppose now that $\al_3-d\xi\geq \gamma$.
Then we have
$$
x^d+a_2\overline{a_1}\equiv 0 \pmod{p^{\al_3-d\xi}},
$$
with $x=x_1'\overline{ x_2'}$, and where $\overline{b}$ denotes the multiplicative
inverse of $b$ modulo $p^{\al_3-d\xi}$. We now appeal to the 
well-known fact that for any $b \in \ZZ$ coprime to $p$, and any $k\geq
\gamma$, 
the number of solutions to the congruence 
$$
x^d\equiv b \pmod{p^k}
$$ 
is either $0$ or $p^{\gamma-\delta-1}\hcf\big(d,p^{\delta}(p-1)\big).
$
This therefore ensures the existence of $K\leq 2^{\gamma-\delta-1}d$ integers
$\la_1,\ldots,\la_{K}$ such that
$$
a_1\la_k^d+a_2\equiv 0 \pmod{p^{\al_3-d\xi}},
$$
for $1\leq k\leq K$. In particular, the point $\x\in\ZZ^3$ in which we are
interested must satisfy $x_1=p^\xi x_1'$, $x_2=p^\xi x_2'$ and 
$$
x_1'\equiv \la_k x_2' \pmod{p^{\al_3-d\xi}},
$$
for some $1\leq k\leq K$. Assuming that $d\geq 2$, these conditions define a union of $K$
lattices, each  of dimension $3$ and determinant
$p^{\al_3+2\xi-d\xi}\geq p^{2\al_3/d}$.
This completes the proof of the claim.
\end{proof}

We are now ready to establish Theorem \ref{upper}. 
Let $a=|a_1a_2a_3|$. 
In view of the lemma, the points that we are
interested in belong to a union of $J\leq 2d^{\omega(a)}$ 
lattices $\sfl_1,\ldots,\sfl_J\subseteq \ZZ^3$, each of dimension $3$
and determinant $\gg a^{2/d}.$ Let us
consider the overall contribution from the vectors belonging to one
such lattice $\sfl_j$, say.  
We will work with a minimal basis $\mathbf{b}^{(1)}, \mathbf{b}^{(2)}, \mathbf{b}^{(3)}$
for $\sfl_j$, which satisfies the well-known bound
$|\mathbf{b}^{(1)}||\mathbf{b}^{(2)}||\mathbf{b}^{(3)}|\gg a^{2/d}$,
and furthermore, whenever 
$\x=\sum_i \la_i\mathbf{b}^{(i)}$ for $\la_i \in \ZZ$,
so it follows that 
$$
\la_i \ll \frac{|\x|}{|\mathbf{b}^{(i)}|}\ll \frac{B}{|\mathbf{b}^{(i)}|}=B_i,
$$
say. Here $|\mathbf{z}|:=\max_{i} |z_i|$ for any $\mathbf{z}\in \RR^3$.
On carrying out this change of variables, \eqref{eq:fermat} becomes $G_j(\la_1,\la_2,\la_3)=0$,
with $G_j$ a ternary form of degree $d$  that is defined over $\ZZ$. We
are now interested in counting integer solutions to this equation with
$\la_i\ll B_i$. It follows from a simple application of Siegel's lemma that
the number of such vectors  is
$$
\ll
1+\Big(\frac{B^3}{|\mathbf{b}^{(1)}||\mathbf{b}^{(2)}||\mathbf{b}^{(3)}|}\Big)^{\frac{1}{2}}
\ll  1+ a^{-\frac{1}{d}}B^{\frac{3}{2}}.
$$
On summing over the $J$ lattices, this therefore establishes Theorem \ref{upper}.

\end{document}